\documentclass[11pt]{article}
\usepackage{epic,latexsym,amssymb}
\usepackage{amsfonts}
\usepackage{amscd}
\usepackage{amsmath}
\usepackage{graphicx}
\usepackage{color}
\usepackage{caption,
caption}
\usepackage{tikz}
\usepackage{amsthm}
\usepackage{bbm}
\usepackage{amsfonts}
\usepackage{xcolor}

\usepackage[margin=0.25in]{geometry}
\usepackage{pgfplots}
\pgfplotsset{width=10cm,compat=1.9}

\usepgfplotslibrary{external}
\usepackage{epic,latexsym,amssymb}
\usepackage{tikz}
\usepackage[bold,full]{complexity}
\usepackage{tikz,epsf}
\usepackage{hyperref}
\usepackage{graphicx}[width=\textwidth]
\usepackage{amsmath}
\usepackage{listings}

\usepackage{amsfonts}
\usepackage{amscd}
\usepackage{amsmath}
\usepackage{graphicx}
\usepackage{color}
\usepackage{caption,subcaption}
\usepackage{lineno}

\usepackage[ruled,vlined]{algorithm2e}

\usepackage{lineno}

\newtheorem{thm}{Theorem}

\newtheorem{prop}[thm]{Proposition}
\newtheorem{lem}[thm]{Lemma}

\newtheorem{cor}[thm]{Corollary}

\newtheorem{claim}{Claim}

\newcommand{\sub}{\mathrm{sub}}
\newcommand{\ecc}{{\rm ecc}}
\newcommand{\diam}{{\rm diam}}
\newcommand{\rad}{{\rm rad}}

\newcommand{\smallqed}{{\tiny ($\Box$)}}

\newcommand{\QEDmark}{\mbox{\textsc{qed}}}
\newcommand{\proofStarter}[1]{\textsc{#1} }

\def\vertex(#1){\put(#1){\circle*{2}}}
\def\vertexo(#1){\put(#1){\circle{2}}}
\def\vert(#1){\put(#1){\circle*{1.5}}}
\def\verto(#1){\put(#1){\circle{1.5}}}
\def\lab(#1)#2{\put(#1){\makebox(0,0)[c]{#2}}}
\definecolor{DarkGreen}{rgb}{0.2, 0.6, 0.3}

\definecolor{electricindigo}{rgb}{0.44, 0.0, 1.0}

\let\oldenumerate\enumerate
\renewcommand{\enumerate}{
  \oldenumerate
  \setlength{\itemsep}{0.5pt}
  \setlength{\parskip}{0pt}
  \setlength{\parsep}{0pt}
}
\begin{document}
\title{Lower bounds for the total (distance) \\ $k$-domination  number of a graph}
\author{$^{1,2}$Randy Davila\
\\
\\
$^1$Research and Development \\
RelationalAI \\
Berkeley, CA 94704, USA\\
\small {\tt Email: randy.davila@relational.ai}\\
\\
$^2$Department of Computational Applied \\ Mathematics \& Operations Research\\
Rice University\\
Houston, TX 77005, USA \\
\small {\tt Email: randy.r.davila@rice.edu} \\
\\
}

\date{}
\maketitle
\begin{abstract}
For $k \geq 1$ and a graph $G$ without isolated vertices, a \emph{total (distance) $k$-dominating set} of $G$ is a set of vertices $S \subseteq V(G)$ such that every vertex in $G$ is within distance $k$ to some vertex of $S$ other than itself. The \emph{total (distance) $k$-domination number} of $G$ is the minimum cardinality of a total $k$-dominating set in $G$, and is denoted by $\gamma_{k}^t(G)$. When $k=1$, the total $k$-domination number reduces to the \emph{total domination number}, written $\gamma_t(G)$; that is, $\gamma_t(G) = \gamma_{1}^t(G)$. This paper shows that several known lower bounds on the total domination number generalize nicely to lower bounds on total (distance) $k$-domination.
\end{abstract}


\section{Introduction}
Distance domination in graphs is a well-known concept, with a recent survey estimating that there are more than 100 papers on the topic to date, including a chapter in the domination monograph by Haynes, Hedetniemi, and Henning~\cite{HaHeHe-topics}. This paper considers a specific, widely applicable variant of distance domination called \emph{total distance $k$-domination}, and from hereon, referred to as \emph{total $k$-domination} for simplicity. This concept was introduced by Henning, Oellermann, and Swart in~\cite{HeOeSw}, where they define a \emph{total $k$-dominating set} of a graph $G$ to be a set $S \subseteq V(G)$ of vertices so that every vertex is within distance $k$ from some vertex of $S$ other than itself. More specifically, let $k \geq 1$ be an integer. A set $S\subseteq V(G)$ is total $k$-dominating in $G$ if for every vertex $v \in V(G)$, we have $d_G(v, S \setminus \{v\}) \leq k$. The \emph{total $k$-domination number} of $G$ is the cardinality of a minimum total $k$-dominating set in $G$ and is denoted by $\gamma_{k}^t(G)$. A total $k$-dominating set of $G$ with cardinality $\gamma_{k}^t(G)$ is called a $\gamma_{k}^t$-set of $G$. Since every vertex $v$ of $G$ must be within distance $k$ to some vertex different from $v$ in a total $k$-dominating set, we note that total $k$-domination is not defined for any graph with isolated vertices. We further remark that when $k = 1$, total $k$-domination reduces to a well-known and heavily studied notion of total domination. In particular, we note that $\gamma_t(G) = \gamma_1^t(G)$, where $\gamma_t(G)$ is the \emph{total domination number} of $G$. For more on total domination, see the excellent monograph by Henning and Yeo~\cite{HeYe-total-domination}.

As with the total domination number of a graph, the computation of the total $k$-domination number is NP-hard~\cite{HaHeHe-topics}. For this reason, many of the results on the total $k$-domination have focused primarily on finding tight upper and lower bounds. For example, the upper bound of Henning, Oellermann, and Swart~\cite{HeOeSw}, which states that $\gamma_{k}^t(G) \leq \frac{2n}{2k+1}$, whenever $G$ is a connected graph with order $n \geq 2k + 1 \geq 3$. However, there do not seem to be many well-known lower bounds for the total $k$-domination number in the literature, which differs from that of the non-total distance dominating variant~\cite{DaFaHeKe}. This observation motivates our contributions. In particular, we extend several known lower bounds on the total domination number to lower bounds on the total $k$-domination number. 

\noindent\textbf{Notation and Terminology.} In this paper, we only consider finite and simple graphs. Let $G$ be a graph with vertex set $V(G)$ and edge set $E(G)$. The order and size of $G$ will be denoted by $n(G) = |V(G)|$ and $m(G) = |E(G)|$, respectively. A nontrivial graph is a graph of order at least 2. Two vertices $v, w \in V(G)$ are neighbors or adjacent whenever $vw \in E(G)$. The open neighborhood of a vertex $v \in V(G)$, written $N_G(v)$, is the set of all neighbors of $v$, whereas the closed neighborhood of $v$ is $N_G[v] = N_G(v) \cup \{v\}$. The degree of a vertex $v \in V(G)$, written $d_G(v)$, is the number of neighbors of $v$ in $G$; and so, $d_G(v) = |N_G(v)|$. A subgraph $H$ of $G$ is a graph where $V(H) \subseteq V(G)$ and $E(H) \subseteq E(G)$. If $H$ is a subgraph of $G$, we write $H \subseteq G$. The complete graph, path, and cycle on $n$ vertices will be denoted by $K_n$, $P_n$, and $C_n$, respectively.

A graph $G$ is connected if, for all vertices $v$ and $w$ in $G$, a $(v, w)$-path exists. A \emph{tree} is a connected graph that contains no cycle as a subgraph. A \emph{forest} is any graph that does not contain a cycle as a subgraph. If a forest is connected, it is necessarily a tree; if not, it is a disjoint union of trees called its components. A vertex of degree 1 in a tree is called a leaf, and a vertex with a leaf neighbor is a \emph{support vertex}. The \emph{distance} from a vertex $v$ to a vertex $w$ in $G$, denoted $d_G(v, w)$, is the length of a shortest $(v, w)$-path in $G$. The distance from a vertex $v$ to a set $S\subseteq V(G)$, denoted $d_G(v, S)$, is the length of a shortest $(v, w)$-path for all $w\in S$. A vertex $v$ is said to \emph{$k$-dominate} a vertex $w$ different than $v$ if the distance from $v$ to $w$ is at most $k$. The \emph{eccentricity} of $v \in G$, written $\ecc_G(v)$ is the distance between $v$ and a vertex farthest from $v$ in $G$. The minimum eccentricity among all vertices of $G$ is the \emph{radius} of $G$, denoted by $\rad(G)$, while the maximum eccentricity among all vertices of $G$ is the \emph{diameter} of $G$, denoted by $\diam(G)$. Thus, the diameter of $G$ is the maximum distance among all pairs of vertices of $G$. A vertex $v$ with $\ecc_G(v) = \diam(G)$ is called a \emph{peripheral vertex} of $G$. A \emph{diametrical path} in $G$ is the shortest path in $G$ whose length is equal to the graph's diameter. Thus, a diametrical path is a path of length $\diam(G)$ joining two peripheral vertices of $G$. 

For notation and graph terminology not introduced here, we refer the reader to~\cite{HaHeHe2024}. We will also use the standard notation $[k] = \{1, \ldots, k\}$.

\section{Preliminaries}
In this section we provide tools fore our main results in the next section. To begin, recall that every total $k$-dominating set of a spanning subgraph of the graph $G$ is a total $k$-dominating set of $G$, and thus, we immediately have the following proposition. 
\begin{prop}\label{prop:basic}
For $k\geq 1$, if $H$ is a spanning subgraph of a graph $G$ without isolated vertices, then 
\[
\gamma_{k}^t(G) \leq \gamma_{k}^t(H).
\]
\end{prop}

We now present our first generalization from total domination to total $k$-domination. More specifically, recall that in 2007 DeLaViña et al.~\cite{DeLiPeWaWe} showed every non-trivial connected graph $G$ has a spanning tree $T$ satisfying $\gamma_t(G) = \gamma_t(T)$. This result generalizes from total domination to total $k$-domination, as shown by the following lemma. 
\begin{lem}\label{lem:spanning-tree}
For $k\geq 1$, if $G$ is a connected graph of order $n\geq 2$, then $G$ has spanning tree $T$ such that 
\[
\gamma_{k}^t(T)= \gamma_{k}^t(G).
\]
\end{lem}
\proof Let $k\geq 1$ and let $G$ be a connected graph of order $n\geq 2$. Next let $S \subseteq V(G)$ be a $\gamma_k^t$-set of $G$. Thus, every vertex $v$ of $G$ is within distance to $k$ to the set $S\setminus \{v\}$ and $\gamma_{k}^t(G) = |S|$. For each $i \in [k]$, define the set $D_i$ to be the set of all vertices $v \in V(G)$ for which $d_G(v, S \setminus \{ v \}) = i$; that is,
\[
D_i = \Big\{ v \in V(G) : d_G(v, S \setminus \{v \}) = i \Big \}.
\]
Since $S$ is a total $k$-dominating set of $G$, every vertex $v \in V(G)$ belongs to some $D_i$ for $i \in [k]$. Moreover, if $v \in D_i$, then $v$ has at least 1 neighbor in $D_{i-1}$, since otherwise $v$ would not be at distance $i$ with $S\setminus \{v \}$. Note that it is also possible that if $v \in D_i$, then $v$ could have neighbors in $D_i$ and $D_{i+1}$. We next construct a spanning subgraph $F \subseteq G$ as follows. For each $i \in [k]$, apply the following rules:
\begin{itemize}
    \item[1.] For each vertex $v \in D_i$, delete all but 1 of the edges joining $v$ to the set $D_{i-1}$. 

    \item[2.] For each vertex $v \in D_i$, delete all edges, if any, that join $v$ to other vertices in the set $D_i$.
\end{itemize}

Since we did not delete any vertices in the above steps, we note that $F$ is a spanning subgraph of $G$. We next show that $F$ is necessarily a spanning forest of $G$. 
\begin{claim}
The subgraph $F$ is a spanning forest.
\end{claim}
\proof By way of contradiction, suppose that $F$ is not a forest. Thus, $F$ contains at least one cycle as a subgraph. Let $C$ be one such cycle in $F$ and let $v$ be a vertex of $C$ of maximum distance (noninclusive) to the set $S$ in $G$; that is, $d_G(v, S\setminus\{v\}) \geq d_G(w, S \setminus\{w\})$ for all $w \in S$. Suppose now that $d_G(v, S \setminus \{ v \}) = \ell$, and so, $v \in D_{\ell}$. If there were a vertex $w$ on $C$ with $d_G(w, S \setminus \{w \}) > \ell$, then $v$ would not be a vertex of $C$ at maximum (noninclusive) distance from $S$, a contradiction. Thus, $d_G(w, S \setminus \{ w \}) \leq \ell$ for all vertices $w$ in $C$. 

Let $v_1$ and $v_2$ be the 2 neighbors of $v$ on the cycle $C$. Next recall that when constructing $F$ from $G$, we removed any edges connecting $v$ to other vertices in $D_{\ell}$. Thus, $v_1$ and $v_2$ cannot be elements of the set $D_{\ell}$. Hence, the neighbors $v_1$ and $v_2$ of $v$ on $C$ are both elements of the set $D_{\ell -1}$, a contradiction since $v$ will have at most one neighbor in the set $D_{\ell -1 }$ by construction of $F$. Therefore, $F$ cannot contain a cycle as a subgraph, and so, $F$ is a spanning forest of $G$. \smallqed

We now construct a spanning tree $T$ from $F$ as follows. If $F$ is a tree, then let $T = F$. Otherwise, $F$ is a disjoint union of $\ell$ tree components for some $\ell \geq 2$. In the case that $F$ consists of $\ell$ components, we let $T$ be the tree obtained from $F$ by adding to $F$ $\ell -1$ so that the resulting graph $T$ is connected. Note that $T$ is necessarily a spanning tree of $G$. Moreover, for $i \in [k]$ and $v \in V(G)$, if $v \in D_i$, then there is a path from $v$ to $S\setminus \{v\}$ of length $i$ in $T$, and so, $d_T(v, S \setminus \{v \}) \leq d_G(v, S \setminus \{v \})$. However, $T$ is a spanning subgraph of $G$, which implies $d_G(v, S \setminus \{v \}) \leq d_T(v, S \setminus \{v \})$ for all $v \in V(G)$. Thus, $T$ is a spanning tree of $G$ that is distance-preserving from $S$, in the sense that $d_T(v, S \setminus \{v \}) = d_G(v, S \setminus \{v \})$ for all $v \in V(G)$. It follows that since $S$ is a total $k$-dominating set of $G$, it must be that $S$ is also a total $k$-dominating set of $T$. Thus, $\gamma_k^t(T) \leq |S| = \gamma_k^t(G)$. However, Proposition~\ref{prop:basic} states $\gamma_{k}^t(G) \leq \gamma_{k}^t(T)$. Hence, $\gamma_{k}^t(G) = \gamma_{k}^t(T)$. \qed

\section{Lower Bounds for $\gamma_k^t(G)$}
In this section we generalize several known lower bounds for the total domination number to those for the total $k$-domination number. To begin, recall that for $k\geq 1$, the \emph{$k$-neighborhood} of a vertex $v \in V(G)$ is the set of all vertices different from $v$ that are at distance at most $k$ to $v$; that is, $N_{G, k}(v) = \{ w : 1 \leq d_G(v, w) \leq k\}$. The \emph{$k$-degree} of a vertex $v \in V(G)$ is defined as $d_{G, k}(v) = |N_{G, k}(v)|$. The \emph{maximum $k$-degree} and \emph{minimum $k$-degree} of a graph $G$ are denoted $\Delta_k(G)$ and $\delta_k(G)$, respectively. The \emph{$k$-degree sequence} of $G$ is the list of all vertex $k$-degree's of $G$ in nonincreasing order and is written $D_k(G) = \Big(\Delta_k(G) = d_k^1, \dots, d_k^n = \delta_k(G) \Big )$. 

For $k\geq1$ and a graph $G$ of order $n$, we define the \emph{sub-total (distance) $k$-domination number} of $G$, written $\sub_k^t(G)$, as the smallest integer $j$ such that $d_k^1 + \dots + d_{k}^j \geq n$; that is,
\[
\sub_k^t(G) = \min\Big\{j : \sum_{i=1}^jd_k^i \geq n \Big \}
\]
When $k = 1$, the sub-total $k$-domination number reduces to the \emph{sub-total domination number}, written $\sub_t(G)$ and introduced by Davila in~\cite{Davila}, and simultaneously introduced by Gentner and Rautenbach in~\cite{GeRa}, where they adopt the notation $sl_t(G)$ to emphasize the connection to the well-known \emph{Slater number} of a graph~\cite{Sl}. Regardless of notation, $\sub_t(G)$ (and $sl_t(G)$) serve as simple lower bounds for the total domination number. That is, $\gamma_t(G) \geq \sub_t(G)$ for any isolate-free graph $G$~\cite{Davila, GeRa}; which is notably a simple application of the \emph{degree sequence index strategy} (\emph{DSI-strategy}) formally introduced by Caro and Pepper in~\cite{CaPe}. The following theorem shows that this bound can easily be generalized from total domination to one for total $k$-domination.
\begin{thm}\label{thm:sub-domination}
For $k \geq 1$, if $G$ is a isolate-free graph, then
\[
\gamma_{k}^t(G) \geq \sub_k^t(G).
\]
\end{thm}
\proof Let $k\geq 1$ and let $G$ be an isolate-free graph of order $n$. Next let $j = \sub_{k}^t(G)$ and let $S \subseteq V(G)$ be a $\gamma_{k}^t$-set of $G$. Thus, every vertex $v$ of $G$ is within distance to $k$ to the set $S\setminus \{v\}$ and $\gamma_{k}^t(G) = |S|$. Hence, $V(G) = \bigcup_{v \in S} N_{G, k}(v)$. Therefore, 
\begin{align*}
n & = |V(G)| \\
  & = \left|\bigcup\limits_{v \in S} N_{G, k}(v) \right| \\
  & \leq \sum\limits_{v \in S} |N_{G, k}(v)| \\
  & = \sum\limits_{v \in S} d_{G, k}(v),
\end{align*}
and so, $\sum_{v\in S}d_{G, k}(v) \geq n$. Since the sum of the $k$-degrees from the vertices in $S$ is at most the sum of the first $|S|$ entries from the $k$-degree sequence $D_k(G)$, we obtain the following inequality,
\[
\sum_{i=1}^{|S|}d_k^i \geq \sum_{v\in S}d_{G, k}(v) \geq n.
\]
Since $j$ is the smallest integer satisfying $\sum_{i = 1}^{j} d_{k}^i \geq n$, it follows that 
\[
\gamma_{k}^t(G) = |S| \geq j = \sub_{k}^t(G),
\]
and the desired inequality is established. \qed 

A trivial lower bound for the total domination number of any isolate-free graph of order $n \geq 2$ and maximum degree $\Delta$ is $\gamma_t(G) \geq \frac{n}{\Delta}$. As a simple application of Theorem~\ref{thm:sub-domination}, we next show that this lower bound on the total domination number also generalizes to a bound for the total $k$-domination number. 
\begin{cor}\label{cor:order-and-degree}
For $k \geq 1$, if $G$ is an isolate-free graph with order $n$ and maximum $k$-degree $\Delta_k$, then
\[
\gamma_{k}^t(G) \geq \frac{n}{\Delta_k}
\]
\end{cor}
\proof Let $G$ be an isolate-free graph of order $n$ and maximum $k$-degree $\Delta_k$. Next let $j = \sub_k^t(G)$, and note that Theorem~\ref{thm:sub-domination} implies $\gamma_{k}^t(G) \geq j$. For each $i \in [n]$ observe that $d_k^i \leq \Delta_k$. Thus, we obtain the following chain of inequalities,
\[
j \Delta_k = \sum_{i = 1}^{j} \Delta_{k} \geq \sum_{i = 1}^{j} d_{k}^i  \geq n. 
\]
Hence, 
\[
\gamma_{k}^t(G) \geq j \geq \frac{n}{\Delta_k},
\]
and the desired bound is established. \qed 

In~\cite{DeLiPeWaWe} DeLaViña et al. proved $\gamma_t(G) \geq \frac{1}{2}(\diam(G)+1)$ for any connected graph $G$ of order $n\geq 2$. The following theorem extends this lower bound from the total domination number to a lower bound on the total $k$-domination number. 
\begin{thm}\label{thm:diam}
For $k \geq 1$, if $G$ is a connected graph of order $n\geq 2$, then 
\[
\gamma_{k}^t(G) \geq \frac{\diam(G)+1}{2k}.
\]
\end{thm}

\proof Let $G$ be a connected graph with order $n\geq 2$ and let $d = \diam(G)$.
Let $P: u_0 u_1 \ldots u_d$ be a diametrical path in $G$, joining two peripheral vertices $u = u_0$ and $v = u_d$ of $G$, and so, $d_G(u, v) = d$ and let $S$ be a $\gamma_{k}^t$-set of $G$ of $G$. Thus, every vertex $v$ of $G$ is within distance to $k$ to the set $S\setminus \{v\}$ and $\gamma_{k}^t(G) = |S|$. We proceed by first proving the following claim about the number of vertices of $P$ that any vertex of $S$ may $k$-dominate. 
\begin{claim}\label{claim:at-most-2k}
Any vertex of $S$ will $k$-dominate at most $2k$ vertices in $P$.
\end{claim}
\proof Let $q \in S$ be a vertex of $S$ that $k$-dominates at least one vertex from the path $P$. Next let $Q$ be the set of all vertices in $P$ that are $k$-dominated by $q$ and then let $i$ and $j$ be the smallest and largest integers in the set $[d]$, respectively, such that $u_i \in Q$ and $u_j \in Q$. Since every vertex of $Q$ is a vertex in the path $P$ and since we choose a smallest integer $i$ with $u_i \in Q$, and largest integer $j$ with $u_j \in Q$, it must be the case that $Q \subseteq \{u_i, u_{i+1}, \ldots, u_j \}$. 

Let $P': u_i \dots d_j$ be the path subgraph of $P$ connecting $u_i$ to $u_j$ in $P$ and note that since $P$ is a shortest $(u, v)$-path in $G$, it must be the case that $P'$ is a shortest $(u_i, u_j)$-path in $G$. Next let $P_i$ be a shortest $(u_i, q)$-path in $G$ and let $P_j$ be a shortest $(q, u_j)$-path in $G$, and observe that since $q$ will $k$-dominate both $u_i$ and $u_j$ in $G$, both paths $P_i$ and $P_j$ have at length most $k$. Therefore, the $(u_i, u_j)$-path, say $P''$, obtained by following the path $P_i$ from $u_i$ to $q$, and then proceeding along the path $P_j$ from $q$ to $u_j$, has length at most $2k$, for otherwise, the vertices $u_i$ and $u_j$ on the path $P$ would not be the vertices with minimum and maximum indices so that $u_i \in Q$ and $u_j \in Q$, respectively. Thus, $P''$ has order at most $2k+1$. Since $q$ does not $k$-dominate itself, $q$ will $k$-dominate at most $(2k+1)-1 = 2k$ vertices on $P''$, which in-turn implies that $q$ will $k$-dominate at most $2k$ vertices on $P'$. Since $Q\subseteq V(P')$, $q$ will $k$-dominate at most $2k$ vertices in the path $P$. \smallqed 

By Claim~\ref{claim:at-most-2k}, each vertex in $S$ will $k$-dominate at most $2k$ vertices of the diametrical path $P$. Further, since $P$ has $d+1$ vertices, we observe $2k|S| \geq  d + 1$.
Hence,
\[
\gamma_{k}^t(G) \geq \frac{d + 1}{2k} = \frac{\diam(G) + 1}{2k},
\]
and the desired bound is established. \qed 

Another bound for the total domination number is $\gamma_t(G) \geq \rad(G)$~\cite{DeLiPeWaWe}. We generalize this radius lower bound for the total domination number to a lower bound for the total $k$-domination number by applying Theorem~\ref{thm:diam} and Lemma~\ref{lem:spanning-tree}.
\begin{thm}\label{thm:radius}
For $k \geq 1$, if $G$ is a connected graph of order $n\geq 2$, then 
\[
\gamma_{k}^t(G) \geq \frac{1}{k}\rad(G).
\] 
\end{thm}
\proof Let $k\geq 1$ and let $G$ be a connected graph of order $n \geq 2$. By Lemma~\ref{lem:spanning-tree}, $G$ has a spanning tree $T$ such that $\gamma_{k}^t(G) = \gamma_{k}^t(T)$. Adding edges to $T$ will not increase the radius, and so $\rad(G)\leq \rad(T)$. Moreover, since $T$ is a tree, $\diam(T) \geq 2\rad(T) - 1$. Thus, by applying Theorem~\ref{thm:diam}, we obtain the following chain of inequalities,
\[
\gamma_{k}^t(G) = \gamma_{k}^t(T) \geq \frac{\diam(T) + 1}{2k} \geq \frac{2\rad(T)}{2k} \geq \frac{1}{k}\rad(G).
\]
Hence, $\gamma_{k}^t(G) \geq \frac{1}{k}\rad(G)$, and the desired bound is established. \qed

\end{document}